\documentclass[10pt]{amsart}

\usepackage{amsfonts}

\usepackage{graphicx}
\usepackage{amsmath}
\usepackage{setspace}
\usepackage{bbm}

\newcommand{\halmos}{{\mbox{\, \vspace{3mm}}} \hfill
\mbox{$\Box$}}

\newtheorem{theorem}{Theorem}

\newtheorem{lemma}[theorem]{Lemma}

\newtheorem{proposition}[theorem]{Proposition}
\newtheorem{remark}[theorem]{Remark}

\def\te#1{\mathrm{e}^{#1}}

\newcommand{\PP}{{\mathbb P}}
\newcommand{\EE}{{\mathbb E}}

\usepackage[english]{babel}
\usepackage[left=3cm,right=3cm,top=2.5cm,bottom=3cm]{geometry}

\usepackage{amsmath,amssymb,amsthm,bbm,color,graphics,version}
\usepackage{mathrsfs}

\newcommand{\bearno}{\begin{eqnarray*}}
\newcommand{\enarno}{\end{eqnarray*}}

\title{Parisian quasi-stationary distributions for asymmetric L\'{e}vy processes}

\author{Irmina Czarna}
\address{Mathematical institute, University of Wroc{\l}aw, Poland.}
\email{czarna@math.uni.wroc.pl}

\author{Zbigniew Palmowski}
\address{Mathematical Institute, University of Wroc{\l}aw, Poland.}
\email{zbigniew.palmowski@gmail.com}

\thanks{This work is partially supported by the Ministry of Science and
Higher Education of Poland under the grants
DEC-2013/09/B/HS4/01496 (2014-2016).
All the authors kindly acknowledge partial support by the project RARE -318984, a Marie Curie IRSES Fellowship within the 7th European Community Framework Programme.}

\date{\today}

\begin{document}

\begin{abstract}
In recent years there has been some focus on quasi-stationary behaviour of an one-dimensional L\'evy process $X$,
where we ask for the law $\PP(X_t\in dy | \tau^-_0>t)$ for $t\to\infty$ and \linebreak
$\tau_0^-=\inf\{t\geq 0: X_t<0\}$.
In this paper we address the same question for so-called Parisian ruin time $\tau^\theta$, that happens when process stays below zero longer than
independent exponential random variable with intensity $\theta$.

\vspace{3mm}


\end{abstract}

\maketitle

\pagestyle{myheadings} \markboth{\sc I.\ Czarna --- Z.\ Palmowski} {\sc Parisian quasi-stationary distributions }

\vspace{1.8cm}

\tableofcontents

\newpage

\section{Introduction}\label{sec:intro}
Let $X=\{X_t : t\geq 0\}$  be a spectrally one-sided L\'evy process defined on the filtered space
$(\Omega, \mathcal{F}, \mathbb{F}, \PP)$ where the filtration $\mathbb{F}=\{\mathcal{F}_t : t\geq 0\}$ is assumed to satisfy
the usual conditions for right continuity and completion.
Suppose now that probabilities $\{\PP_{x}\}_{x \in \mathbb{R}}$ corresponds to the
conditional version of $\PP$ where $X_{0}=x$ is given. We simply write $%
\PP_{0}=\PP $. % We assume that $X_t\to-\infty$ a.s. as $t\to\infty$.

Define the first passage time into the lower half line $(-\infty,0)$ by
\[
\tau_0^-=\inf\{t \geq 0: X_t <0\}.
\]

In recent years there has been some focus on the existence and characterization of the so-called
 limiting quasi-stationary distribution (or Yaglom's limit) defined by:
\begin{equation}\label{def2}\mu(dy):=\lim_{t\uparrow\infty}\PP_x(X_t\in dy|\tau_0^- >t).\end{equation}
The sense in which this limit is quasi-stationary follows the classical interpretations of works such as Seneta and Vere-Jones \cite{Seneta},
Tweedie \cite{tweedie}, Iglehart \cite{Iglehart} (for a random walk), Jacka and Roberts \cite{jackaroberts}, Kyprianou \cite{kyprianou} (within the context of the $M/G/1$ queue), Martinez and San Martin \cite{Martmart} (for a Brownian motion with drift),  Kyprianou and Palmowski \cite{KP} (for a general light-tailed L\'{e}vy process), Hass and Rivero \cite{rivero} (for regularly varying L\'{e}vy process), Mandjes et al. \cite{MPR} (for a workload process of single server queue) and other references therein.

In this paper, the principal object of interest is the quasi-stationary distribution of Parisian type, where
$\tau^-_0$ is replaced by so-called Parisian ruin time:
\begin{equation}\label{definition:parisiantime}\tau^{\theta}=\inf\{t > 0: t-g_t > \textbf{e}_{\theta}^{g_t}, X_t<0 \},\end{equation}
where $g_t=\sup\{s\leq t: X_s \geq 0\}$ and $\textbf{e}_{\theta}^{g_t}$ is an independent of $X$ exponential random variable with intensity $\theta>0$ related to a separate negative excursion $g_t$.
The ruin time $\tau^\theta$ happens when process $X_t$ stays negative longer than $\textbf{e}_{\theta}^{g_t}$, which we will refer as implementation clock.
We want to emphasize that in the definition of $\tau^{\theta}$ there is not a single underlying random variable but a whole sequence of independent copies of a generic exponential random variable $\textbf{e}_{\theta}^{g_t}$ each one of them attached to a separate excursion below zero. The model with exponentially distributed delay has also been studied by Landriault et al. \cite{LJFZ} and by Baurdoux et al. \cite{BPPR}. The name for this ruin comes from Parisian option that prices
are activated or canceled depending on type of option if
underlying asset stays above or below barrier long enough in a row (see \cite{Hansjoerg, chesneyetal, DassWu3}).
So far only probability of Parisian ruin is known (see \cite{DassWu2} and \cite{czarnapalmowski}).
In this paper we will find sufficient conditions for existence and identify
%via its Laplace transform
the following limit:
\begin{equation}\label{qspar}\mu^\theta_x(dy):=\lim_{t\uparrow\infty}\PP_x(X_t\in dy|\tau^\theta >t).\end{equation}

%In particular we prove that for a spectrally positive L\'evy process:
%\[\mu(dy)=...,\]
%where ...
%Similarly, for a spectrally negative L\'evy process we have
%\[\mu(dy)=...,\]
%where ...
%Note that some expression derived in this paper also identify (after inverting respective Laplace transform)
%the finite-time Parisian ruin probability.

The idea of the proof of the main results is based on finding double Laplace transform of \linebreak $\PP_x(X_t\in dy, \tau^\theta >t)$ with respect to space and time.
Then for some specific form of the L\'evy measure (that will be defined later) using 'Heavyside' operation we will identify the asymptotics of this probability as $t\to\infty$ (see e.g. \cite{AW1997} and \cite{henrici}).
%This allows us to identify the quasi-stationary distribution (\ref{qspar}).

The paper is organized as follows.
In the next section we state some preliminary facts. Later in Theorem \ref{masterThm} we give the formula, which is essential to obtain the main results of this paper. In Section \ref{sec:main} we present the quasi-stationary distribution for one-sided L\'evy processes.
% and analyze some special examples.
Finally, in the Appendix the proof of Theorem \ref{masterThm} is presented.

\section{Preliminaries}\label{sec:prel}
\subsection{Asymmetric L\'evy processes}
In this section we present definitions, notations and basic facts on L\'evy processes
(we also refer to Bertoin \cite{bertoin1996} and Kyprianou \cite{Kbook} for a complete introduction to the theory of L\'evy processes).
We will focus on asymmetric L\'evy processes, which are either spectrally positive (having nonnegative jumps) or spectrally negative (having nonpositive jumps). We danote by $\Pi_X(\cdot)$ the L\'evy measure of $X$.
With $X$ we associate the Laplace exponent $\varphi(\beta):=\frac{1}{t}\log \EE(e^{\beta X_t})$
defined for all $\beta$ for which exists, where
$$
\varphi(\beta) = \gamma \beta + \frac{1}{2} \sigma^2 \beta^2 + \int_{\mathbb{R}} \left( \mathrm{e}^{-\beta z} - 1 + \beta z \mathbbm{1}_{(-1,1)}(z) \right) \Pi_X(\mathrm{d}z) ,
$$
for $\gamma \in \mathbb{R}$ and $\sigma \geq 0$.
Moreover, we define function $\Phi(q)=\sup\{\beta \geq 0: \varphi(\beta)=q\}$ called right-inverse of $\varphi$.
We will also consider so-called dual process $\widehat{X}_t=-X_t$ with the L\'evy measure $\Pi_{\widehat{X}}\left(0,y\right)=\Pi_{X}\left(-y,0\right)$.
Characteristics of $\widehat{X}$ will be indicated by
using a hat over the existing notation for characteristics of $X$.
In particular, the probabilities $\widehat{\PP}_x$ and the expectations $\widehat{\EE}_x$ concern the dual process.

For the process $X$ we define the ascending ladder height process $(L^{-1},H)=\{(L^{-1}_t,H_t)\}_{t\geq 0}$:
$$L^{-1}_t:=
\left\{ {\inf\{s>0: L_s>t \} \quad \textrm{ if $t < L_{\infty}$} \atop \infty \quad \textrm{otherwise}} \right.
$$
and
$$H_t:= \left \{
X_{L_t^{-1}} \quad \textrm{if $t < L_{\infty}$} \atop \infty \quad \textrm{otherwise},
\right.$$
where $L\equiv\{L_t\}_{t\geq 0}$ is the local time at the maximum (see \cite[p. 140]{Kbook}).
Recall that $(L^{-1}_t,H_t) $ is a (killed) bivariate subordinator with the Laplace exponent
$\kappa(\alpha,\beta)=-\frac{1}{t}\log \EE\left (e^{-\alpha L_t^{-1}-\beta H_t}\mathbbm{1}_{\{t\leq L_{\infty}\}}\right)$ and with the jump measure
$\Pi_H$.
We define the descending ladder height process $(\widehat{L}^{-1},\widehat{H})=\{(\widehat{L}^{-1}_t,\widehat{H}_t)\}_{t\geq 0}$
with the Laplace exponent $\widehat{\kappa}(\alpha, \beta)$ constructed from dual process $\widehat{X}$.
%Recall that $L_\infty$ has exponential distribution with parameter $\kappa(0,0)$.
Moreover, for a spectrally negative L\'evy process from the Wiener-Hopf factorization we have:
\begin{equation}\label{kappy}
\kappa(\alpha, \beta)=\Phi(\alpha)+\beta,\qquad \widehat{\kappa}(\alpha, \beta)=\frac{\alpha-\varphi(\beta)}{\Phi(\alpha)-\beta};
\end{equation}
see \cite[p. 169-170]{Kbook}.
%Hence
%\begin{equation}\label{kappazero}
%\widehat{\kappa}(0,0)=\varphi^\prime(0+).
%\end{equation}

We introduce the
%potential measure $\mathcal{U}$ defined by
%$$ \mathcal{U}(dx,ds)=\int_0^{\infty} \P(L^{-1}_t \in ds, H_t \in dx) dt$$
%with the Laplace transform
%$ \int_{[0,\infty)^2} e^{-\theta s -\beta x} \mathcal{U}(dx,ds)=1/\kappa(\theta,\beta) $ and
renewal function:
$$V(dx)=\int_{[0,\infty)}\mathcal{U}(dx,ds)=\EE\left(\int_0^{\infty} \mathbbm{1}_{\{H_t\in dx\}}\, dt \right).$$
In particular,
\begin{equation}
\int_0^\infty e^{-\alpha z} V (z) \,dz =\frac{1}{\alpha\kappa(0,\alpha)}.
\end{equation}
For spectrally negative L\'{e}vy process upward ladder height process is
a linear drift and hence the renewal measure is just the Lebesgue measure:
\begin{equation}\label{renfunleb}
V(dx)=dx.
\end{equation}
Moreover, from \cite[p. 195]{Kbook}:
\begin{equation}\label{LTdlahatU}
\int_0^\infty e^{-\alpha z}\widehat{V}(dz)=\frac{\alpha}{\varphi(\alpha)}.
\end{equation}

We will also use the first passage times:
$$\tau_x^-=\inf\{t\geq 0: X_t < x\},\qquad \tau_x^+=\inf\{t\geq 0: X_t \geq x\}.$$
 We define Girsanov-type change of measure via:
\begin{equation}\label{Girsanov}
\left. \frac{d\PP_{x}^{c}}{d\PP_{x}}\right| _{\mathcal{F}_{t}}=\frac{\mathcal{E}%
_{t}\left( c\right) }{\mathcal{E}_{0}\left( c\right) }
\end{equation}
for any $c$ for which $\EE e^{cX_1}<\infty$, where $\mathcal{E}_{t}\left( c\right) =\exp \{cX_{t}-\varphi
\left( c\right) t\}$ is exponential martingale under $\PP_x$ and $\mathcal{F}_{t}$ is a natural filtration of $X$.
It is easy to check that under this new
measure $X$ still remains within the class of L\'evy processes with the
Laplace exponent
\begin{equation}\label{varphinowe}
\varphi _{c}\left( \beta \right) :=\varphi \left( \beta +c\right) -\varphi \left(
c\right).
\end{equation}
%for $\theta \geq -c.$
All quantities calculated for $\PP^{c}$ will have subindex $c$ added to their counterparts
defined on $\PP$.
%By (\ref{kappy}) on $\PP^c$ exponent of descending ladder process equals,
%\begin{equation}\label{kappynowe}
%\widehat{\kappa}_c(0,\beta)=\frac{\varphi \left( \beta +c\right) -\varphi \left(
%c\right)}{\beta}
%\end{equation}
%since on an account of positive drift of $X$ we have $\Phi (0)=0$.
\subsection{Scale functions and Exit problems}
For $q\geq 0$, there exists a function $W^{(q)}: [0,\infty) \to [0,\infty)$,
called {\it $q$-scale function}, that is continuous and
increasing with the Laplace transform:
\begin{equation}\label{eq:defW} \int_0^\infty
\te{-\alpha x} W^{(q)} (x) \, d x = \frac{1}{\varphi(\alpha) - q},\quad\alpha >
\Phi(q). \end{equation}
We denote $W^{(0)}(x)=W(x)$.
Domain of $W^{(q)}$ is extended to the entire real
axis by setting $W^{(q)}(x)=0$ for $x<0$. For each $x\geq 0$, function $q\to W^{(q)}(x)$ may be analytically extended to $q\in \mathcal{C}$.
Moreover, let
$$Z^{(q)}(x)=1+q\int_0^x W^{(q)}(y)\,dy.$$
The initial value of $W^{(q)\prime}$ for spectrally negative process is given by
\begin{equation}\label{initialvalues}
W^{(q)\prime}(0+) =
\begin{cases}
2/\sigma^2 & \text{when $\sigma>0$,} \\
\left(\Pi(0,\infty)+q\right)/\left(\gamma+\int_{(0,1)}z\Pi_X(dz)\right)^2 & \text{when $\sigma=0$ and $\Pi(0,\infty)<\infty$,} \\
\infty & \text{otherwise.}
\end{cases}
\end{equation}
Finally for spectrally negative process we recall the following identity taken from
% $\cite{Kbook} and
\cite{LJF}:
\begin{eqnarray}%\label{exit}
%\EE_x\left[e^{-q\tau_y^+}, \tau_y^+<\infty\right]&=&e^{-\Phi(q)(y-x)},\\
%\EE_{x}\left[ e^{-q\tau _{0}^{-}}, \tau_{0}^{-}<\infty\right] &=&Z^{(q)}(x)-\frac{q}{\Phi \left( q\right) }%
%W^{(q)}(x)\;,  \label{one-sided-down}\\
\label{exit2} \\\nonumber
\EE_x\left[e^{-q\tau _{0}^{-}+vX_{\tau^-_0}},\tau_0^- <\infty\right]
&=&e^{vx}\left(1 +\left(q-\varphi(v)\right)\int_{0}^{x} e^{- vy} W^{(q)}(y) dy \right)-\frac{q-\varphi(v)}{\Phi(q)-v}W^{(q)}(x).
\end{eqnarray}
Note that when $v=\Phi(q)$ the therm $\frac{q-\varphi(v)}{\Phi(q)-v}$ should be read as
$\lim_{v\rightarrow \Phi(q)}\frac{q-\varphi(v)}{\Phi(q)-v}=1/\Phi'(q)$.

\subsection{Classical results}
Kyprianou and Palmowski \cite{KP} and  Mandjes et al. \cite{MPR} proved that the quasi-stationary distribution (\ref{def2}) related with the classical ruin time $\tau_0^-$ equals
%\begin{equation*}
%\mu(dy) = \xi^{\star} \kappa_{\xi^{\star}} (0,\xi^{\star})e^{- \xi^{\star} y}V_{\xi^{\star}}(y)dy \;\textrm{ on $[0,\infty)$}.
%\end{equation*}
\begin{itemize}
\item if $X$ is a spectrally negative process
\[
\mu(dy) = (q^{\star})^2 y e^{-q^{\star} y} \, dy,
\]
where $\varphi(\cdot)$ attains its strictly negative minimum at $q^{\star}>0$ and hence $\varphi'(q^{\star})=0$.
\item if $X$ is a spectrally positive process
\[
\mu(dy) = Q_+  e^{q^{\star} y} V_{-q^{\star}}(y) \,dy,
\]
where $Q_+:=\left(\int_0^\infty e^{-q^{\star} z} V_{-q^{\star}} (z) \, dz\right)^{-1}$
and $\widehat{\varphi}(\cdot)$ attains its strictly negative minimum at $q^{\star}<0$ and hence $\widehat{\varphi}'(q^{\star})=0$.
\end{itemize}
The key lemma in deriving above results comes from the Wiener-Hopf factorization (see \cite[Eqn. (10)]{KP} and \cite[Thm 4.]{MPR}) and holds true for any spectrally one-sided L\'evy process:
\begin{lemma}\label{KeyFact}
For $x\geq 0$, $\alpha \in \mathbb{R}$
\begin{eqnarray*}
\EE_x\left[e^{-\alpha X_{\textbf{e}_q}},\tau^-_0>\textbf{e}_q\right]=\int_0^{\infty}q e^{- q t} \EE_x\left[e^{-\alpha X_{t}},\tau^-_0>t\right] dt=
\frac{\kappa(q,0)}{\kappa(q,\alpha)}e^{-\alpha x}\left(\int_{[0,x]} e^{\alpha z } \widehat{\PP}\left(\widehat{\overline{X}}(e_{q}) \in dz\right)\right),
\end{eqnarray*}
where $\widehat{\overline{X}}(t) =\sup_{s\leq t} \widehat{X}(s)$ and $\textbf{e}_q$ is an exponential random variable with intensity $q$, independent of X. 
\end{lemma}
In particular for spectrally positive L\'{e}vy processes
\begin{equation}\label{E_classic_SP}
\EE_x\left[e^{-\alpha X_{\textbf{e}_q}},\tau^-_0>\textbf{e}_q\right]=\frac{q}{q-\widehat{\varphi}(\alpha)}\left(e^{-\alpha x}-e^{-\widehat{\Phi}(q)x}\right)
\end{equation}
and for spectrally negative processes
\begin{equation}\label{E_classic_SN}
\begin{split}
\EE_x\left[e^{-\alpha X_{\textbf{e}_q}},\tau^-_0>\textbf{e}_q\right]=&\frac{q}{\Phi(q)+\alpha}\left(W^{(q)}(x)- e^{-\alpha x}W^{(q)}(0)+W^{(q)}(0)\right)-q e^{-\alpha x}
\int_0^{x} e^{\alpha z} W^{(q)}(z)\,dz.
\end{split}
\end{equation}

\subsection{Tauberian-type results}
To obtain main result we will specialize Theorem \cite[Th. 37.1]{Doetsch1974}
in the following way. For $\mathfrak W$-contour
with an half-angle of opening
$\pi/2<\psi\leq \pi$ (defined formally e.g. in \cite[Fig 30, p. 240]{Doetsch1974}) we will consider a cone
${\mathcal G}_{\alpha}(\psi)$ which is a region between contour $\mathfrak W$
and line $\Re(z)=0$.

\begin{proposition}(Heaviside principle). \label{t.tauberian}
 Suppose that $\tilde f(z)=\int_0^\infty e^{-z x} f(x) \;dx$ (defined for some function $f$) satisfies the following three conditions for some $\alpha < 0$:
\begin{itemize}\item [\bf{(A1)}] $\tilde f(z)$ is analytic in the
region ${\mathcal G}_{\alpha}(\psi)$,
\item [\bf{(A2)}] $\tilde f(z) \to 0$ as $|z| \to \infty$ for $z \in {\mathcal G}_{\alpha}(\psi)$,
\item [\bf{(A3)}]  for some constant $K$ and non integer real number $s$
\begin{eqnarray*}
\label{eqn:bridge case 1}
\tilde f(z)=K\;\mathbbm{1}(s > 0) - C (z-\alpha)^{s}+o((z-\alpha)^{s}),
\end{eqnarray*}
for ${\mathcal G}_{\alpha}(\psi)\ni z \to \alpha$.
\end{itemize}
  Then
\begin{eqnarray*}
f(x)=\frac{C}{\Gamma(-s)} x^{-s-1}e^{\alpha x}(1+o(1))\;,
\end{eqnarray*}
  where $K$ must be $\tilde f(\alpha)$ if $s > 0$ and recall that some function $h(x) = o(1)$ if $\lim_{x\rightarrow \infty} h(x)=0.$
\end{proposition}

\section{Intermediate results}
In this section we present formulas for the double Laplace transforms, which are essential to obtain the main results.
\begin{theorem}\label{masterThm}
For $\alpha, x \geq 0$,
\begin{itemize}
\item if $X$ is a spectrally positive L\'evy process then
\begin{equation} \label{Spectrpositive}
\begin{split}
\EE_x\left[e^{-\alpha X_{\textbf{e}_q}},\tau^\theta>\textbf{e}_q\right]
= & e^{-\widehat{\Phi}(q) x}\frac{q^2\left(\widehat{\Phi}(\theta+q)-\widehat{\Phi}(q)\right)}{\theta\left(q+\theta\right)\left(\widehat{\Phi}(q)+\alpha  \right)}\\&+
\frac{q}{q-\widehat{\varphi}(\alpha)}\left(e^{-\alpha x}-e^{-\widehat{\Phi}(q) x}\frac{\left(q+\theta-\widehat{\varphi}(\alpha)\right)\left( \widehat{\Phi}(\theta+q)-\widehat{\Phi}(q)\right)}{\theta\left(\widehat{\Phi}(\theta+q)-\alpha\right)}\right),
\end{split}
\end{equation}
where $\textbf{e}_q$ is an exponential random variable with intensity $q$, independent of X and  $\textbf{e}_{\theta}^{g_t}$. 
\item if $X$ is a spectrally negative L\'evy process then
\begin{equation}\label{Spectrnegative}
\begin{split}
\EE_x\left[e^{-\alpha X_{\textbf{e}_q}},\tau^\theta>\textbf{e}_q\right]=&\EE_x\left[e^{-\alpha X_{\textbf{e}_q}},\tau^-_0>\textbf{e}_q\right]
+\frac{q^2}{\left(q-\widehat{\varphi}(\alpha)\right)(q+\theta)}\EE_x\left[e^{-q\tau_0^- -\alpha X_{\tau^-_0}},\tau_0^-<\infty \right] \\
&-\frac{q^2}{\left(q-\widehat{\varphi}(\alpha) \right)(q+\theta)}\left(e^{\Phi(q)x}-\frac{1}{\Phi'(q)}\right)\\
&+\EE\left[e^{-\alpha X_{\textbf{e}_q}},\tau^\theta>\textbf{e}_q\right]\EE_x\left[e^{-q\tau_0^- +\Phi(\theta+q)X_{\tau^-_0}},\tau_0^-<\infty \right],
\end{split}
\end{equation}
where \begin{equation}\label{Ezeroboundedn}
\EE\left[e^{-\alpha X_{\textbf{e}_q}},\tau^\theta>\textbf{e}_q\right]
=\frac{q (\Phi(\theta+q)-\Phi(q))}{(q+\theta)(\Phi(q)+\alpha)\theta}+\frac{q^2 (\Phi(\theta+q)-\Phi(q))}{(q+\theta)\Phi'(q)(q-\widehat{\varphi}(\alpha))\theta}
\end{equation} and appropriate expected values are given in (\ref{exit2}) and (\ref{E_classic_SN}).
\end{itemize}
\end{theorem}
The proof of Theorem \ref{masterThm} is given in Appendix.
\section{Main results}\label{sec:main}
In this section we use the Laplace transforms given in Theorem \ref{masterThm} to identify the quasi-stationary distribution for the spectrally one-sided cases.
\subsection{Spectrally positive case}
In this section we assume that $X$ is a spectrally positive L\'{e}vy process satisfying the following conditions:\\
There exists $q_+<0$ such that
\begin{itemize}
\item[\bf{(SP1)}] $\widehat{\varphi}(q)<\infty$ for $ q_+ <q$;
\item[\bf{(SP2)}] $\widehat{\varphi}(q)$ attains its negative minimum at $q^\star <0$, where $q_+<q^\star<0$ and hence $\widehat{\varphi}'(q^\star)=0$;
%\item[(SP3)] $X_1$ is non-lattice.
\item[\bf{(SP3)}] Function $\widehat{\Phi}$ can be analytically extended into ${\mathcal G}_{\xi^*}(\psi)$ for some \linebreak $\pi /2<\psi\leq\pi $ and $\xi^{\star}:=\widehat{\varphi}(q^\star)<0$.
\end{itemize}

\begin{remark} \rm
Since $\widehat\Phi(\vartheta)$ is the Laplace exponent of a subordinator
we have the following spectral
representation:
\begin{equation}\label{reprhatphi}\widehat\Phi({\vartheta})=d_+{\vartheta}+
\int_0^\infty (1-e^{-{\vartheta} y})\,\Pi_+(dy)\end{equation}
with $\int_0^\infty (y\wedge 1)\Pi_+(dy)<\infty$. From its definition we see that $\xi^\star$
must be a singular point of $\widehat{\Phi}$. Moreover, if
there exists a density of $\Pi_+$ which is of semiexponential type, 
then condition $\bf{(SP3)}$ is satisfied. 
We recall that function $f$  is said to be {\em semiexponential} if for some  $0< \phi\le\pi/2$, there exists finite
and strictly negative
$\gamma({\vartheta})$, defined  as the infimum of all such $a$ such that
\[\left|f(e^{i{\vartheta}}r)\right|<e^{ar}\]
for all sufficiently large $r$ and
$-\phi\le {\vartheta}\le \phi$.
In particular, this assumption holds for example
for a linear Brownian motion $X(t)=\sigma B(t)-c t$, where $c>0$ (see for details \cite{MPR}).
\end{remark}
Before we present the main formula of this section, we first recall a few important facts.
\begin{lemma}\label{expansion}
Under $\bf{(SP1)-(SP3)}$,
\begin{equation}\label{Phi_MPR1}
\widehat{\Phi}(q)=q^{\star}+k^{\star}(q-\xi^{\star})^{1/2}+o((q-\xi^{\star})^{1/2}),
\end{equation}
as $q\downarrow \xi^{\star}$, where $k^\star:=\sqrt{{2}/{ \widehat{\varphi}^{''}(q^\star)}}$.
\proof For details see \cite[Lem. 10]{MPR}. \endproof
\end{lemma}

\begin{proposition}\label{qasmain}
As $q\downarrow \xi^{\star}$,
\[\EE_x\left[e^{-\alpha X_{\textbf{e}_q}},\tau^\theta>\textbf{e}_q\right]= C_p(\alpha, x)+H_p(\alpha, x)(q-\xi^*)^{1/2}+{\rm o} ((q-\xi^*)^{1/2})\]
for some $C_p$ and
\begin{equation}
\begin{split}
H_p(\alpha, x)=&\frac{A_p(x)}{\left(\alpha-\widehat{\Phi}(\theta+\xi^{\star}) \right)\left(\widehat{\varphi}(\alpha)-\xi^{\star} \right)}
- \frac{A_p(x) e^{-q^\star x}}{\theta\left(\alpha-\widehat{\Phi}(\theta+\xi^{\star}) \right)}\\
+&B_p(x) \left(-x\frac{\widehat{\Phi}(\theta+\xi^{\star})-q^\star}{q^\star+\alpha}-\frac{\widehat{\Phi}(\theta+\xi^{\star})}{(q^\star+\alpha)^2}\right),
\end{split}
\end{equation}
where $A_p(x)=\xi^{\star}\left(\widehat{\Phi}(\theta+\xi^{\star})xk^\star+k^\star-x k^\star\right)$ and
$B_p(x)=\frac{e^{-q^\star x} (\xi^{\star})^2}{\theta(\theta+\xi^{\star})}$.
\proof The above result is obtained using Theorem \ref{masterThm} (Eqn. \eqref{Spectrpositive}) together with Lemma \ref{expansion} and simple facts saying
that if functions $a$ and $b$ have expansions
$a(q)=a^{\star}+a^\circ (q-\xi^{\star})^{1/2}+o((q-\xi^{\star})^{1/2})$ and
$b(q)=b^{\star}+b^\circ (q-\xi^{\star})^{1/2}+o((q-\xi^{\star})^{1/2})$
then
\begin{eqnarray}
a(q)b(q)&=&a^{\star}b^{\star}+(a^\circ b^\star+b^\circ a^\star) (q-\xi^{\star})^{1/2}+o((q-\xi^{\star})^{1/2}),\label{raz}\\
\frac{a(q)}{b(q)}&=&\frac{a^*}{b^*} + \left(\frac{a^\circ}{b^\star} -\frac{b^\circ}{(b^\star)^2} a^\star\right) (q-\xi^{\star})^{1/2}+o((q-\xi^{\star})^{1/2}),\label{dwa}\\
e^{a(q)x}&=& e^{a^\star x} +a^\circ x e^{a^\star x} (q-\xi^{\star})^{1/2}+o((q-\xi^{\star})^{1/2}).\label{trzy}
\end{eqnarray} \endproof
\end{proposition}

Hence by Proposition \ref{qasmain} and Proposition \ref{t.tauberian} we obtain that
\begin{equation}\label{Hp}
\lim_{t\rightarrow \infty}\EE_x[e^{-\alpha X_t}, \tau^-_0>t]=
H_p(\alpha,x) \frac{t^{-3/2}}{\Gamma(-1/2)}e^{\xi^\star t}(1+o(1))
\end{equation}
 and by setting $\alpha=0$ we have
\begin{equation}\label{Hp0}
\lim_{t\rightarrow \infty}\PP_x( \tau^-_0>t)=
H_p(0,x) \frac{t^{-3/2}}{\Gamma(-1/2)}e^{\xi^\star t}(1+o(1)).
\end{equation}

The main result of this section is the following theorem.
\begin{theorem}\label{mainspecpos}
Assume that conditions $\bf{(SP1)-(SP3)}$ hold. Then the Parisian quasi-stationary distribution~(\ref{qspar}) equals
\begin{equation*}
\begin{split}
\mu_x^{\theta}(dy)=& \frac{1}{H_p(0,x)}\left(A_p(x)\left(e^{\widehat{\Phi}(\theta+\xi^{\star})y}\ast \widehat{W}^{(\xi^\star)}(y)\right)
-\frac{A_p(x)e^{-q^\star x}}{\theta}e^{\widehat{\Phi}(\theta+\xi^{\star})y}\right)\\
-&\frac{B_p(x)e^{-q^\star x}}{H_p(0,x)}\left(\widehat{\Phi}(\theta+\xi^{\star})(x+y)-x q^\star \right),
\end{split}
\end{equation*}
where $\ast$ denotes convolution, $H_p(0,x)$, $A_p(x)$ and $B_p(x)$ are given in Proposition \ref{qasmain}.
\proof
From Equations \eqref{Hp} and \eqref{Hp0} we get that the Parisian quasi-stationary distribution (\ref{qspar}) has the Laplace transform $H_p(\alpha,x)/H_p(0,x)$. The statement of the theorem is obtained by inverting above Laplace transform with respect to $\alpha$ (for details see \cite{Abramowitz}). \endproof
\end{theorem}

\subsection{Spectrally negative case}
In this section we assume that $X$ is a spectrally negative L\'{e}vy process satisfying the following conditions:\\
There exists $q_->0$ such that
\begin{itemize}
\item[\bf{(SN1)}] $\varphi(q)<\infty$ for $0< q <q_-$;
\item[\bf{(SN2)}] $\varphi(q)$ attains its negative minimum at $q^\star >0$, where $0<q^\star<q_-$ and hence $\varphi'(q^\star)=0$;
%\item[(SP3)] $X_1$ is non-lattice.
\item[\bf{(SN3)}] Function $\widehat{\Phi}$ can be analytically extended into ${\mathcal G}_{\xi^*}(\psi)$ for some \linebreak $\pi /2<\psi\leq\pi $ and $\xi^{\star}:=\varphi(q^\star)<0$.
\end{itemize}

\begin{remark} \rm
Since $\Phi(\vartheta)$ is the Laplace exponent of a subordinator
we have the following spectral
representation:
\begin{equation}\label{reprhatphiSN}\Phi({\vartheta})=d_-{\vartheta}+
\int_0^\infty (1-e^{-{\vartheta} y})\,\Pi_-(dy)\end{equation}
and $\int_0^\infty (y\wedge 1)\Pi_-(dy)<\infty$. From its definition we see that $\xi^\star$
must be a singular point of $\Phi$. Moreover, if
there exists a density of $\Pi_-$ which is of semiexponential type, then
condition $\bf{(SN3)}$ is satisfied.
\end{remark}
As for spectrally positive processes, before we present the main formula of this section, we first present a few important facts.
\begin{lemma}\label{expansionSN}
Under $\bf{(SN1)-(SN3)},$
\begin{equation}\label{Phi_MPR2}
\Phi(q)=q^{\star}+k^{\star}(q-\xi^{\star})^{1/2}+o((q-\xi^{\star})^{1/2}),
\end{equation}
as $q\downarrow \xi^{\star}$, where $k^\star:=\sqrt{{2}/{ \varphi^{''}(q^\star)}}$.
\proof For details see \cite[Lem. 16]{MPR}. \endproof
\end{lemma}

\begin{proposition}\label{qasmainSN}
As $q\downarrow \xi^{\star}$
\[\EE_x\left[e^{-\alpha X_{\textbf{e}_q}},\tau^\theta>\textbf{e}_q\right]= C_n(\alpha, x)+H_n(\alpha, x)(q-\xi^*)^{1/2}+{\rm o} ((q-\xi^*)^{1/2})\]
for some $C_n$ and
\begin{equation}
\begin{split}
H_n(\alpha, x)=&\frac{A_n(x)}{(q^\star+\alpha)^2}+B_n(x)\left(\frac{1}{q^\star+\alpha}+\frac{\Phi(\theta+\xi^\star)- q^\star}{(q^\star+\alpha)^2}\right)\\
+&\frac{k^\star W^{(\xi^{\star})}(0)e^{-\alpha x}}{(q^\star+\alpha)^2}
-\frac{(\xi^\star)^2 e^{-q^\star x}k^\star x }{(\xi^\star+\alpha)(\varphi(-\alpha)-\xi^\star)},
\end{split}
\end{equation}
where $A_n(x)=-\xi^{\star}k^{\star}\left(W^{(\xi^{\star})}(x)+W^{(\xi^{\star})}(0) +
\frac{\xi^{\star}-1}{\xi^{\star}+\theta}W^{(\xi^{\star})}(x)\right)$ and \\
$B_n(x)=\frac{-\xi^{\star}k^{\star}e^{\Phi(\theta+\xi^{\star})x}}{\xi^{\star}+\theta}\left(1-\theta \int_0^x e^{\Phi(\theta+\xi^{\star})y} W^{(\xi^{\star})}(y) dy \right)$.
\proof The above result is obtained using Theorem \ref{masterThm} (Eqn. \eqref{Spectrnegative}) together with Lemma \ref{expansionSN} and 
facts \eqref{raz}-\eqref{trzy}.
\endproof
\end{proposition}

Hence by Proposition \ref{qasmainSN} and Proposition \ref{t.tauberian} we obtain that
\begin{equation}\label{HpSN}
\lim_{t\rightarrow \infty}\EE_x[e^{-\alpha X_t}, \tau^-_0>t]=
H_n(\alpha,x) \frac{t^{-3/2}}{\Gamma(-1/2)}e^{\xi^\star t}(1+o(1))
\end{equation}
 and by setting $\alpha=0$ we have
\begin{equation}\label{Hp0SN}
\lim_{t\rightarrow \infty}\PP_x( \tau^-_0>t)=
H_n(0,x) \frac{t^{-3/2}}{\Gamma(-1/2)}e^{\xi^\star t}(1+o(1)).
\end{equation}

The main result of this section is the following theorem.
\begin{theorem}\label{mainspecposSN}
Assume that conditions $\bf{(SN1)-(SN3)}$ hold. Then the Parisian quasi-stationary distribution~(\ref{qspar}) equals
\begin{equation*}
\begin{split}
\mu_x^{\theta}(dy)=&\frac{e^{-q^\star y}}{H_n(0,x)}\left(A_n(x) y +B_n(x)+\left(\Phi(\theta+\xi^\star)-q^\star \right)y\right) \\+ & \frac{k^\star}{H_n(0,x)}\left( W^{(\xi^\star)}(0) (y-x) u(y-x)e^{-q^\star(y- x)}-(\xi^{\star})^2e^{-q^\star x} \left(e^{\xi^{\star}y}\ast W^{(-\xi^\star)}(y)\right)\right),
\end{split}
\end{equation*}
where $\ast$ denotes convolution, $u(\cdot)$ is the Heaviside step function, $A_n(x)$ and $B_n(x)$ are given in Proposition \ref{qasmainSN}.
\proof
From Equations \eqref{HpSN} and \eqref{Hp0SN} we get that the Parisian quasi-stationary distribution (\ref{qspar}) has the Laplace transform $H_n(\alpha,x)/H_n(0,x)$. The statement of the theorem is obtained by inverting above Laplace transform with respect to $\alpha$ (for details see \cite{Abramowitz}). \endproof
\end{theorem}

\section*{Appendix}
In this section we present the proof of the formula given in Theorem \ref{masterThm}.
First, we assume that the process X is spectrally positive i.e. we will prove Equation \eqref{Spectrpositive}.

{\it Proof of Theorem \ref{masterThm}, Eq. \eqref{Spectrpositive}. }

For $\epsilon \geq 0 $ and $g_t=\sup\{s\leq t: X_s \geq 0\}$ we define the stopping time:
$$\tau_{-\epsilon}^{\theta}=\inf\{t>0: X_{t}-g_t >e^{g_t}_{\theta}, X_{t-e^{g_t}_{\theta}}<-\epsilon \},$$
which is the first time when an excursion starting when $X$ is below $-\epsilon$ ending at the moment of getting back up over zero and of length greater than $e^{g_t}_{\theta}$, has occurred. Clearly, we have $\tau_{0}^{\theta}=\tau^{\theta}$. 

Using strong Markov property for $x\geq 0$ we obtain:
\begin{eqnarray}\label{markovpositive}
\EE_x\left[e^{-\alpha X_{\textbf{e}_q}},\tau_{-\epsilon}^\theta>\textbf{e}_q\right]=
\EE_{x+\epsilon}\left[e^{-\alpha X_{\textbf{e}_q}},\tau^-_{0}>\textbf{e}_q\right]
+\PP_{x+\epsilon}(\tau_{0}^-<\textbf{e}_q)\EE_{-\epsilon}\left[e^{-\alpha X_{\textbf{e}_q}}, \tau_{-\epsilon}^\theta >\textbf{e}_q\right].
\end{eqnarray}

Straightforward from \eqref{E_classic_SP} we have:
\begin{equation}\label{markovpositiveclassic}
\EE_{x+\epsilon}\left[e^{-\alpha X_{\textbf{e}_q}},\tau^-_0>\textbf{e}_q\right]
= \frac{q}{q-\widehat{\varphi}(\alpha)}\left(e^{-\alpha (x+\epsilon)}-e^{-\widehat{\Phi}(q) (x+\epsilon)}\right)
\end{equation}
and for $\alpha=0$
\begin{equation}\label{markovpositiveclassic_alpha=0}
\PP_{x+\epsilon}(\tau_{0}^-<\textbf{e}_q)= e^{-\widehat{\Phi}(q) (x+\epsilon)}.
\end{equation}

Moreover, to identify  $\EE_{-\epsilon}\left[e^{-\alpha X_{\textbf{e}_q}},\tau_{-\epsilon}^\theta>\textbf{e}_q\right]$ note that:
\begin{eqnarray*}
\lefteqn{\EE_{-\epsilon}\left[e^{-\alpha X_{\textbf{e}_q}},\tau_{-\epsilon}^\theta>\textbf{e}_q\right]}\\&&= \PP(e_{\theta}>\textbf{e}_q)
\EE_{-\epsilon}\left[e^{-\alpha X_{\textbf{e}_q}},\tau_0^+>\textbf{e}_q\right] +
\int_0^\infty\PP_{-\epsilon}(\tau_0^+\leq \min( \textbf{e}_q,e_{\theta}),  X_{\tau^+_0}\in dz)
\EE_z\left[e^{-\alpha X_{\textbf{e}_q}},\tau_{-\epsilon}^\theta>\textbf{e}_q\right]
\\&&=\frac{q}{q+\theta}
\widehat{\EE}_{\epsilon}\left[e^{\alpha \widehat{X}_{\textbf{e}_q}},\widehat{\tau}_0^->\textbf{e}_q\right] +
\int_0^\infty\PP_{-\epsilon}(\tau_0^+\leq  e_{q+\theta},  X_{\tau^+_0}\in dz)
\EE_z\left[e^{-\alpha X_{\textbf{e}_q}},\tau_{-\epsilon}^\theta>\textbf{e}_q\right].
\end{eqnarray*}
Now using (\ref{markovpositive}) together with \eqref{markovpositiveclassic} and \eqref{markovpositiveclassic_alpha=0} we obtain:
\begin{eqnarray}\label{expe}
\\\nonumber
\lefteqn{\EE_{-\epsilon}\left[e^{-\alpha X_{\textbf{e}_q}},\tau_{-\epsilon}^\theta>\textbf{e}_q\right]}
\\\nonumber&&=
\frac{q}{q+\theta}
\widehat{\EE}_{\epsilon}\left[e^{\alpha \widehat{X}_{\textbf{e}_q}},\widehat{\tau}_0^->\textbf{e}_q\right] \\\nonumber&&+
\frac{q}{q-\widehat{\varphi}(\alpha)}\left\{
e^{-\alpha \epsilon}\widehat{\EE}_{\epsilon}\left[e^{-(q+\theta)\widehat{\tau}_0^- +\alpha \widehat{X}_{\widehat{\tau}^-_0}},\widehat{\tau}_0^-<\infty \right]
- e^{-\widehat{\Phi}(q) \epsilon}\widehat{\EE}_{\epsilon}\left[e^{-(q+\theta)\widehat{\tau}_0^- +\widehat{\Phi}(q) \widehat{X}_{\widehat{\tau}^-_0}},\widehat{\tau}_0^-<\infty \right] \right\}\\\nonumber&&+
e^{-\widehat{\Phi}(q) \epsilon}\EE_{-\epsilon}\left[e^{-\alpha X_{\textbf{e}_q}},\tau_{-\epsilon}^\theta>\textbf{e}_q\right]
\widehat{\EE}_{\epsilon}\left[e^{-(q+\theta)\widehat{\tau}_0^- +\widehat{\Phi}(q) \widehat{X}_{\widehat{\tau}^-_0}},\widehat{\tau}_0^-<\infty \right],
\end{eqnarray}
where $\widehat{\EE}_{\epsilon}\left[e^{\alpha \widehat{X}_{\textbf{e}_q}},\widehat{\tau}_0^->\textbf{e}_q\right]$ is given in
\eqref{E_classic_SN}.

Finally we split the analysis into two cases, when the process $X$ is of bounded variation (BV) and
unbounded variation (UBV). In the first scenario we take $\epsilon=0$ in (\ref{expe}) and since
$\widehat{W}(0)>0$ from \eqref{exit2} we obtain that
\begin{eqnarray}\label{Eeps}
\\\nonumber
\EE\left[e^{-\alpha X_{\textbf{e}_q}},\tau^\theta>\textbf{e}_q\right]
=\frac{q}{q-\widehat{\varphi}(\alpha)}\left(1-\frac{\left(q+\theta-\widehat{\varphi}(\alpha)\right)\left( \widehat{\Phi}(\theta+q)-\widehat{\Phi}(q)\right)}{\theta\left(\widehat{\Phi}(\theta+q)-\alpha\right)}\right)
+\frac{q^2\left(\widehat{\Phi}(\theta+q)-\widehat{\Phi}(q)\right)}{\theta\left(q+\theta\right)\left(\widehat{\Phi}(q)+\alpha  \right)}.
\end{eqnarray}
This Equation together with \eqref{markovpositive} for $\epsilon=0$ proves \eqref{Spectrpositive} in BV case.
For the process of unbounded variation paths (or equivalently when $\widehat{W}(0)=0$) we will use identities \eqref{exit2} and \eqref{E_classic_SN} and the following inequalities:
$$0\leq \lim_{\epsilon \downarrow 0} \frac{\int_0^{\epsilon} e^{-v y}\widehat{W}^{(q+\theta)}(y) dy}{\widehat{W}^{(q+\theta)}(\epsilon) }\leq \lim_{\epsilon \downarrow 0} \int_0^{\epsilon} e^{-v y}dy=0$$
and
$$0\leq \lim_{\epsilon \downarrow 0} \frac{\int_0^{\epsilon} e^{\alpha y}\widehat{W}^{(q)}(y) dy}{\widehat{W}^{(q+\theta)}(\epsilon) }\leq \lim_{\epsilon \downarrow 0} \frac{\widehat{W}^{(q)}(\epsilon)}{\widehat{W}^{(q+\theta)}(\epsilon)}\int_0^{\epsilon} e^{\alpha y}dy=0$$
that follow from l'H\^opital's rule, the fact that $\widehat{W}$ is increasing and the equality \eqref{initialvalues}.
Now taking the limit $\epsilon \downarrow 0$ on both sides of (\ref{expe}) 
will produce Equation \eqref{Eeps} which completes the proof. One can also see that, for $0<\epsilon'<\epsilon$,
$$
\{\tau_{-\epsilon'}^{\theta}>\textbf{e}_q\} \subset \{\tau_{-\epsilon}^{\theta}>\textbf{e}_q\} \quad \text{and} \quad \bigcap_{\epsilon>0}\{\tau_{-\epsilon'}^{\theta}>\textbf{e}_q\}=\{\tau^{\theta}>\textbf{e}_q\} ,
$$
so
\begin{equation*}
\EE_x\left[e^{-\alpha X_{\textbf{e}_q}},\tau^\theta>\textbf{e}_q\right]=\lim_{\epsilon\downarrow0} \EE_x\left[e^{-\alpha X_{\textbf{e}_q}},\tau_{-\epsilon}^\theta>\textbf{e}_q\right] .%=\frac{\mathbb{E} \left[ X_1 \right]}{g_a(r)} .
\end{equation*}

\halmos

{\it Proof of Theorem \ref{masterThm}, Eq. \eqref{Spectrnegative}.}\\
The proof is based on strong Markov property and fact that spectrally negative L\'evy process creeps upward. Note that:
\begin{eqnarray}\label{glowny0}
\lefteqn{\EE_x\left[e^{-\alpha X_{\textbf{e}_q}},\tau^\theta>\textbf{e}_q\right]}\\\nonumber&&=\EE_x\left[e^{-\alpha X_{\textbf{e}_q}},\tau^-_0>\textbf{e}_q\right]\nonumber\\
&&\quad+ \int_0^\infty\PP_x(\tau_0^-\leq \textbf{e}_q, -X_{\tau^-_0}\in dz)\nonumber\\\nonumber&&
\quad\quad\left\{\EE\left[e^{-\alpha X_{\textbf{e}_q}},\tau_z^+ > \textbf{e}_q \right]\PP(e_\theta>\textbf{e}_q)\right.\\\nonumber&&\left.
\quad\quad+\PP(\tau_z^+\leq \min( e_\theta, \textbf{e}_q))\EE\left[e^{-\alpha X_{\textbf{e}_q}},\tau^\theta>\textbf{e}_q\right]\right\}.
\label{glowny}\end{eqnarray}
In (\ref{glowny0}) we use lack of memory of the exponential distribution.

Then
\begin{eqnarray*}
\lefteqn{\EE_x\left[e^{-\alpha X_{\textbf{e}_q}},\tau^\theta>\textbf{e}_q\right]
%&=&\EE_x\left[e^{-\alpha X_{\textbf{e}_q}},\tau^-_0>\textbf{e}_q\right]\nonumber\\
%&&+\int_0^\infty\PP_x(\tau_0^-\leq \textbf{e}_q, -X_{\tau^-_0}\in dz)\nonumber\\&&
%\left\{\frac{q}{q+\theta}\EE\left[e^{-\alpha X_{\textbf{e}_q}},\tau_z^+ > \textbf{e}_q \right]\nonumber\right.\\&&\left.
%+e^{-\Phi(\theta+q)z}\EE\left[e^{-\alpha X_{\textbf{e}_q}},\tau^\theta>\textbf{e}_q\right]\right\} \nonumber
=\EE_x\left[e^
{-\alpha X_{\textbf{e}_q}},\tau^-_0>\textbf{e}_q\right]}\nonumber\\
&&+\frac{q}{q+\theta}\int_0^\infty \PP_x(\tau_0^-\leq \textbf{e}_q, -X_{\tau^-_0}\in dz)\EE_{-z}\left[e^{-\alpha X_{\textbf{e}_q}},\tau_0^+ \geq \textbf{e}_q \right]\nonumber\\&&
+\EE\left[e^{-\alpha X_{\textbf{e}_q}},\tau^\theta>\textbf{e}_q\right]\int_0^\infty e^{-\Phi(\theta+q)z} \PP_x(\tau_0^-\leq \textbf{e}_q, -X_{\tau^-_0}\in dz). \nonumber\\
\end{eqnarray*}
Thus
\begin{eqnarray*}
\lefteqn{\EE_x\left[e^{-\alpha X_{\textbf{e}_q}},\tau^\theta>\textbf{e}_q\right]=\EE_x\left[e^{-\alpha X_{\textbf{e}_q}},\tau^-_0>\textbf{e}_q\right]}\nonumber\\
&&+\frac{q}{q+\theta}\int_0^\infty\PP_x(\tau_0^-\leq \textbf{e}_q, -X_{\tau^-_0}\in dz)\widehat{\EE}_z\left[e^{\alpha \widehat{X}_{\textbf{e}_q}},\widehat{\tau}_0^- \geq \textbf{e}_q \right]\nonumber\\&&
+\EE\left[e^{-\alpha X_{\textbf{e}_q}},\tau^\theta>\textbf{e}_q\right]\EE_x\left[e^{-q\tau_0^- +\Phi(\theta+q)X_{\tau^-_0}},\tau_0^-<\infty \right] \nonumber.
\end{eqnarray*}
In the next step we use (\ref{E_classic_SP}) for $\widehat{\EE}_z\left[e^{\alpha \widehat{X}_{\textbf{e}_q}},\widehat{\tau}_0^- \geq \textbf{e}_q \right]$ to derive the following expression:
\begin{eqnarray}\label{exp_intermediate}
\lefteqn{\EE_x\left[e^{-\alpha X_{\textbf{e}_q}},\tau^\theta>\textbf{e}_q\right]=\EE_x\left[e^{-\alpha X_{\textbf{e}_q}},\tau^-_0>\textbf{e}_q\right]}\\\nonumber
&&+\frac{q^2}{\left(q-\varphi(-\alpha)\right)(q+\theta)}\EE_x\left[e^{-q\tau_0^- -\alpha X_{\tau^-_0}},\tau_0^-<\infty \right] \\\nonumber\\\nonumber&&
-\frac{q^2}{\left(q-\varphi(-\alpha) \right)(q+\theta)}\EE_x\left[e^{-q\tau_0^- +\Phi(q)X_{\tau^-_0}},\tau_0^-<\infty \right]\\\nonumber\\\nonumber&&
+\EE\left[e^{-\alpha X_{\textbf{e}_q}},\tau^\theta>\textbf{e}_q\right]\EE_x\left[e^{-q\tau_0^- +\Phi(\theta+q)X_{\tau^-_0}},\tau_0^-<\infty \right], \nonumber
\end{eqnarray}
where appropriate expected values are given in \eqref{E_classic_SN} and \eqref{exit2} respectively.
Finally, to identify $\EE\left[e^{-\alpha X_{\textbf{e}_q}},\tau^\theta>\textbf{e}_q\right]$
we use the same arguments like for the spectrally positive processes. We split the analysis into two cases, when the process $X$ is of bounded variation (BV) and unbounded variation (UBV). For BV processes we take $x=0$ in (\ref{exp_intermediate}) and we get the result. For UBV processes, for $\epsilon \geq 0$ and $g_t^{\epsilon}=\sup\{0\leq s\leq t:X_s\geq\epsilon\}$ we introduce the stopping time 
$\tau^\theta_\epsilon=\inf\{ t>0: t-g_t^{\epsilon}>e^{g_t^{\epsilon}}_\theta, X_t<0  \}$ which approximates Parisian ruin time $\tau^\theta$, when $\epsilon \downarrow 0$ and finally from \eqref{exp_intermediate} with $x=\epsilon$ and $\tau_{\epsilon}^\theta$ by taking $\epsilon \downarrow 0$ we get \eqref{Ezeroboundedn} and then \eqref{Spectrnegative}.
%The details are left to the reader.

\halmos

%\section*{Acknowledgements}
%This work is partially supported by the Ministry of Science and
%Higher Education of Poland under the grants
%DEC-2013/09/B/HS4/01496 (2014-2016).
%All the authors kindly acknowledge partial support by the project RARE -318984, a Marie Curie IRSES Fellowship within the
%7th European Community Framework Programme.

\end{document}